\documentclass[12pt,twoside]{article}
\usepackage{amsmath}
\usepackage{amsfonts}
\usepackage{amsthm}
\usepackage{amssymb}
\usepackage{graphicx}
\usepackage{multicol}

\begin{document}
\title{{\normalsize{\bf Classical 
Poincar{\'e} conjecture via 4D topology}}}
\author{{\footnotesize Akio KAWAUCHI}\\
{\footnotesize{\it Osaka Central Advanced Mathematical Institute, 
Osaka Metropolitan University}}\\
{\footnotesize{\it Sugimoto, Sumiyoshi-ku, Osaka 558-8585, Japan}}\\
{\footnotesize{\it kawauchi@omu.ac.jp}}}
\date\, 
\maketitle
\vspace{0.25in}
\baselineskip=10pt
\newtheorem{Theorem}{Theorem}[section]
\newtheorem{Conjecture}[Theorem]{Conjecture}
\newtheorem{Lemma}[Theorem]{Lemma}
\newtheorem{Sublemma}[Theorem]{Sublemma}
\newtheorem{Proposition}[Theorem]{Proposition}
\newtheorem{Corollary}[Theorem]{Corollary}
\newtheorem{Claim}[Theorem]{Claim}
\newtheorem{Definition}[Theorem]{Definition}
\newtheorem{Example}[Theorem]{Example}

\begin{abstract} 
The classical Poincar{\'e} conjecture that every homotopy 
3-sphere is diffeomorphic to the 3-sphere is confirmed by Perelman in arXiv papers 
solving Thurston's program on geometrizations of 3-manifolds. A new confirmation of this conjecture is given by a method of 4D topology. For this proof, 
the spun torus-knot of every knot in every homotopy 3-sphere is observed to be a ribbon torus-knot in the 4-sphere, where 
Smooth 4D Poincar{\'e} Conjecture and Ribbonness of a sphere-link with (not necessarily meridian-based) free fundamental group are used. 
By examining a disk-chord system of a ribbon solid torus bounded by the spun torus-knot, 
it is proved that the knot belongs to a 3-ball in the homotopy 3-sphere. Then 
by Bing's result, it is confirmed that the homotopy 3-sphere is diffeomorphic to the 
3-sphere. 

\phantom{x}

\noindent{\footnotesize{\it Keywords:} Homotopy 3-sphere,\, Spun torus-knot,
Ribbon solid torus-knot.}

\noindent{\footnotesize{\it Mathematics Subject Classification 2010}:
Primary 57M40; Secondary 57N13, 57Q45}

\end{abstract}

\baselineskip=15pt

\bigskip

\noindent{\bf 1. Introduction}

A {\it homotopy 3-sphere} is a smooth 3-manifold $M$ homotopy equivalent to 
the 3-sphere $S^3$. The following Poincar{\'e} Conjecture \cite{Poi1,Poi2} is positively shown by Perelman in arXiv papers \cite{P1,P2} solving positively Thurston's program \cite{T} on geometrizations of 3-manifolds (see \cite{M} for detailed historical notes).

\phantom{x} 

\noindent{\bf Poincar{\'e} Conjecture.}  Every homotopy 3-sphere $M$ is diffeomorphic to $S^3$.

\phantom{x} 

A new confirmation of this result is presented here by combining Smooth 4D Poincar{\'e} Conjecture and Free Ribbon Lemma for an $S^2$-link in the 4-sphere $S^4$ with R. H. Bing's result \cite{Bing1,Bing2} on Poincar{\'e} Conjecture.
A  {\it homotopy 4-sphere} is a smooth 4-manifold $X$  homotopy equivalent to the 4-sphere $S^4$.  The following conjecture was a folklore conjecture. 
	
\phantom{x} 

\noindent{\bf Smooth 4D Poincar{\'e} Conjecture.} 
Every smooth homotopy 4-sphere $X$ is diffeomorphic to $S^4$. 

\phantom{x} 
	
The positive proof of this conjecture is shown in \cite{K10}. 
A {\it surface-link} in $S^4$ is a surface $L$ smoothly embedded in $S^4$. When  $L$ is connected, it is a {\it surface-knot}. If all components of 
$L$ are 2-spheres, then it is an $S^2$-{\it link}. 
 A surface-link $L$ in $S^4$ is {\it trivial} if $L$ bounds disjoint handlebodies in $S^4$, and  
 a {\it ribbon surface-link} if $L$ is equivalent to a surface-link obtained  from a trivial $S^2$-link $O$ by surgery along disjointedly embedded 1-handles 
 on $O$ in $S^4$. The following lemma is shown in \cite{K11} as {\it Free Ribbon Lemma} and used in Section~3. 

\phantom{x}

\noindent{\bf Free Ribbon Lemma.} 
Any $S^2$-link $L$ in $S^4$ with free fundamental group $\pi_1(S^4\setminus L,b)$ is a ribbon $S^2$-link in $S^4$.
 
\phantom{x}

The proof of this lemma is moved from this preprint version to the paper \cite{K11} (for completeness of the argument), which is done by  
 using Smooth 4D Poincar{\'e} Conjecture and Smooth Unknotting Conjecture explained as follows: 

\phantom{x}

\noindent{\bf Smooth Unknotting Conjecture.} 
Every smooth surface-link $L$ in $S^4$ with a meridian-based free fundamental group 
$\pi_1(S^4 - L,b)$ is a trivial surface-link. 

\phantom{x}

The proof of this conjecture is shown by \cite{K7,K8,K9}.  
Artin's spinning construction of a knot $k$ in $S^3$ in \cite{Artin} to construct 
the spun $S^2$-knot $K(k)$ in the 4-sphere $S^4$ allows us to generalize to a connected graph $\gamma$ in every homotopy 
3-sphere $M$ to construct the spun $S^2$-link $K(\gamma)$ in a homotopy 4-sphere $X(M)$ which is diffeomorphic to $S^4$ by Smooth 4D Poincar{\'e} Conjecture, so that $X(M)$ is identified with $S^4$. This construction is applied to a Heegaard graph $\gamma$ of $M$ 
(associated to a Heegaard splitting of $M$). 
Then the spun $S^2$-link $K(\gamma)$ is an
$S^2$-link in $X(M)$ with free fundamental group (not always meridian-based free group). 
By Free Ribbon Lemma,  the spun $S^2$-link $K(\gamma)$ is a ribbon $S^2$-link in $X(M)$. 
It is observed that for every knot $k$ in every homotopy 3-sphere $M$, there is  a Heegaard graph $\gamma$ of $M$ 
such that $k$ is contained in the loop system of $\ell(\gamma)$ of $\gamma$. This means that  
the spun $S^2$-knot $K(k)$ of every knot $k$ in every homotopy 3-sphere $M$ is a ribbon $S^2$-knot in $X(M)$. 
Then, by definition, 
the spun torus-knot $T(k)$ of every knot $k$ in every homotopy 3-sphere $M$ is a ribbon torus-knot in $X(M)$. 
Thus, the spun torus-knot $T(k)$ always bounds a ribbon solid torus $V_R$ in $X(M)$. 
By an argument of a disk-chord system of $V_R$ bounded by the spun torus-knot $T(k)$  in $X(M)$, 
the following result is shown.

\phantom{x}

\noindent{\bf Theorem~1.1.} Every knot $k$ in every homotopy 3-sphere $M$ belongs to a 3-ball $D^3$ in $M$.

\phantom{x} 

By combining Theorem~1.1 with the following result of Bing in \cite{Bing1,Bing2}, 
it is proved that every homotopy 3-sphere $M$ is diffeomorphic to $S^3$. 
Thus, the proof of Poincar{\'e} conjecture is completed.

\phantom{x} 

\noindent{\bf Bing's Theorem.} 
A homotopy 3-sphere $M$ is diffeomorphic to $S^3$ if every knot $k$ in $M$ 
belongs to a 3-ball in $M$.

\phantom{x} 

Outline of  the proof of  Poincar{\'e} Conjecture is as follows:

\medskip

\noindent{\bf (1st Step)} By using Smooth 4D Poincar{\'e} Conjecture, show that Artin's spinning construction of every Heegaard graph $\gamma$ of every homotopy 3-sphere $M$ gives a spun $S^2$-link $K(\gamma)$  in $S^4$ with free fundamental group (not always meridian-based free group).

\medskip

\noindent{\bf (2nd Step)} By Free Ribbon Lemma, the spun $S^2$-link $K(\gamma)$ is a ribbon $S^2$-link in $S^4$.

\medskip

\noindent{\bf (3rd Step)} Show that every knot $k$ in $M$ is contained in a loop system $\ell(\gamma)$ of a Heegaard graph $\gamma$ of $M$, so that the 
spun $S^2$-knot $K(k)$ of $k$  is a ribbon $S^2$-knot in $S^4$.
 
\medskip 
 
\noindent{\bf (4th Step)} By definition of a ribbon surface-knot, show that the spun torus-knot $T(k)$ of $k$ in $M$ is a ribbon torus-knot in $S^4$. 

\medskip

\noindent{\bf (5th Step)}  By using a ribbon solid torus $V_R$ bounded by the spun torus-knot  $T(k)$ in $S^4$ and  a disk-chord system of $V_R$, show that $K$ 
belongs to a 3-ball  $D^3$ in $M$.

\medskip

\noindent{\bf (6th Step)} By Bing's theorem, $M$ is diffeomorphic to $S^3$.

\phantom{x}

In Section~2, Artin's spinning construction of a connected graph in a 
homotopy 3-sphere is explained. In Section~3, an argument of a disk-chord system 
of a ribbon solid torus bounded by a ribbon torus-knotis  explained. 
In Section~4, the proof of Theorem~1.1  is done.

\phantom{x}

\noindent{\bf 2. Artin's spinning construction of a connected graph in a 
homotopy 3-sphere}

Throughout this section, $M$ denotes a homotopy 3-sphere unless otherwise mentioned. 
For a homotopy 3-sphere $M$, let $M^{(o)}$ be 
the compact once-punctured manifold $\mbox{cl}(M\setminus B)$ of $M$ 
for a 3-ball $B$ in $M$. 
Let 
\[S=\partial B=\partial M^{(o)}\]
be the boundary 2-sphere of $M^{(o)}$. 
The closed smooth 4-manifold $X(M)$ defined by 
\[X(M)=M^{(o)}\times S^1\cup S\times D^2\]
is called the {\it spun manifold} of $M$ with {\it axis} 4-submanifold 
$S\times D^2$.
As a convention, the 3-submanifold $M^{(o)}\times 1$ of the product 
$M^{(o)}\times S^1$ is identified with $M^{(o)}$. In particular, 
a point $(q,1)\in M^{(o)}\times 1$ is identified with the point $q\in M^{(o)}$. 
This 4-manifold $X(M)$ is a smooth homotopy 4-sphere by the van Kampen theorem and a homological argument 
and hence $X(M)$ is diffeomorphic to the 4-sphere $S^4$ by Smooth 4D Poincar{\'e} Conjecture. {\it From now on, the identification} $X(M)=S^4$ {\it is fixed.}
A {\it legged loop } with {\it base point} $v$ is 
the union $k\cup\omega$ of a loop $k$ and an arc $\omega$ joining 
the base point $v$ with a point of $k$. The arc $\omega$ is called a {\it leg}. 
A {\it legged loop system} with base point $v$ is the union
\[\gamma=\cup_{i=1}^n \ell_i\cup \omega_i\]
of $n$ legged loops $\ell_i\cup \omega_i\, (i=1,2,\dots,n)$ meeting only 
at the same base point $v$. 
Let $\ell(\gamma)=\cup_{i=1}^n \ell_i=\ell_*$ denote the loop system of the legged loop system 
$\gamma$. Let $\omega_*=\cup_{i=1}^n \omega_i$ and 
$v_*=\ell_*\cap \omega_*$. 
A regular neighborhood $B$ of $\omega_*$ in $M$ is taken as a 
3-ball $B$ used for the compact once-punctured manifold $M^{(o)}=\mbox{cl}(M\setminus B)$ of $M$. 
Deform the subgraph $\gamma\cap B$ of $\gamma$ so that 
\[\omega_*\subset B, \quad \omega_*\cap S=v_*
\quad\mbox{and}\quad \ell_*\cap B=\ell_*\cap S=a'_*\]
for a regular neighborhood arc system $a'_*$ of $v_*$ in $\ell_*$. Let 
\[a(\gamma)=\cup_{i=1}^n a_i=a_*\] 
for a proper arc $a_i=\mbox{cl}(\ell_i\setminus a'_i)\,(i=1,2,\dots,n)$ in $M^{(o)}$. 
Let 
\[\dot a(\gamma)=\partial a_*=\partial a'_*\]
be the set of $2n$ points in the boundary 2-sphere $S$ of $M^{(o)}$. 
The {\it spun} $S^2$-link of the graph $\gamma$ 
is the $S^2$-link $K(\gamma)$ in the 4-sphere $X(M)$ defined by 
\[K(\gamma)=a(\gamma)\times S^1\cup\dot a(\gamma)\times D^2.\] 

\phantom{x}

\noindent{\bf Lemma~2.1.} The inclusion 
$M^{(o)}\setminus a(\gamma) \subset X(M)\setminus K(\gamma)$ 
induces an isomorphism
\[\sigma:\pi_1(M\setminus \gamma, v^+)\to \pi_1(X(M)\setminus K(\gamma), v^+)\] 
sending a meridian system of the proper arc system $a(\gamma)$ in $M^{(o)}$ 
to a meridian system of $K(\gamma)$, where the base point $v^+$ is taken in $S\setminus a_*$

\phantom{x}

\noindent{\bf Proof of Lemma~2.1.} 
Note that there is a canonical isomorphism 
\[\pi_1(M^{(o)}\setminus a(\gamma),v^+)\cong\pi_1(M\setminus \gamma,v^+).\] 
Then the desired isomorphism $\sigma$ is obtained by applying the van Kampen theorem between 
$(M^{(o)}\setminus a(\gamma))\times S^1$ 
and $(S\setminus \dot a(\gamma))\times D^2$. 
This completes the proof of Lemma~2.1.
$\square$

\phantom{x}

Here is a note on Lemma~2.1. 

\phantom{x}

\noindent{\bf Note~2.2.} A general connected graph $\gamma$ with Euler characteristic 
$\chi(\gamma)=1-n$ in $M$ 
is deformed into a legged loop system $\gamma$ in $M$ by choosing a maximal tree to shrink to a base point $v$. 
Note that there are only finitely many maximal trees of $\gamma$ such 
that the loop systems $\ell(\gamma)$ of the resulting legged loop systems $\gamma$ 
are distinct as links. By Lemma~2.1, we can obtain finitely 
many distinct spun $S^2$-links in $S^4$ with isomorphic fundamental groups obtained 
by taking different maximal trees of the connected graph $\gamma$. This is a detailed explanation on the spun $S^2$-link of a connected graph 
associated with a maximal tree in \cite[p.204]{K1} when $M=S^3$. 

\phantom{x}

When a homotopy 3-sphere $M$ is given by a Heegaard spitting $V\cup V'$ pasting along a Heegaard surface $F=\partial V=\partial V'$ of genus $n$, a  legged loop system $\gamma$ with  loop system $\ell(\gamma)$ of $2n$ loops is constructed 
as follows. 
A {\it spine} of a handlebody $V$ of genus $n$ is a legged loop system $\gamma_V$ 
in $F=\partial V$ with base point $v$ such that the inclusion map 
$\gamma_V\to V$ induces an isomorphism $\pi_1(\gamma,v)\to \pi_1(V,v)$. 
A regular neighborhood $\dot V$ of $\gamma_V$ in $F$ is a planar surface in $F$. 
By \cite[Theorem~10.2]{Hempel}, 
there is a diffeomorphism $ (\dot V\times [0,1], \dot V\times 0)\to (V,\dot V)$ 
sending every point $(x,0)\in\dot V\times 0$ to $x\in\dot V$.
The surface $\dot V$ is called a {\it spine surface} of $V$. 
Let $\gamma_V$ and $\gamma_{V'}$ be spines of the handlebodies $V$ and $V'$ in 
$F$ with the same base point $v$, respectively. 
A {\it Heegaard graph} of $M$ is a legged loop system 
$\gamma=\gamma_M$ in $M$ with base point $v$ which is the union of legged loop systems $\gamma^+_V$ and $\gamma^+_{V'}$
obtained from $\gamma_V$ and $\gamma_{V'}$ by pushing 
$\gamma_V\setminus v$ and $\gamma_{V'}\setminus v $ into the interiors 
$\mbox{Int}V$ and $\mbox{Int}V'$, respectively. 
The following lemma is obtained. 

\phantom{x}

\noindent{\bf Lemma~2.3.} For every Heegaard graph $\gamma$ of every homotopy 
3-sphere $M$, the fundamental group 
$\pi_1(X(M)\setminus K(\gamma),v^+)$ of the spun $S^2$-link $K(\gamma)$ in the 
4-sphere $X(M)$ is a free group of rank $2n$. 

\phantom{x}

\noindent{\bf Proof of Lemma~2.3.} The closed complement 
$\mbox{cl}(M\setminus N(\gamma))$ for a regular neighborhood $N(\gamma)$ 
of $\gamma$ in $M$ is diffeomorphic to the handlebody 
$F^{(o)}\times[-1,1]$ for the once-punctured 
surface $F^{(o)}$ of $F$. Since the fundamental group 
$\pi_1(F^{(o)}\times[0,1],v^+)$ with base point $v^+$ taken in $(\partial F^{(o)})\times[0,1]$ is a free group of rank $2n$, 
the desired result is obtained from Lemma~2.1.
$\square$

\phantom{x}

It is noted that this free group in Lemma~2.3 is not necessarily a meridian-based free group. Here is an example. 

\phantom{x}

\begin{figure}[hbtp]
\begin{center}
\includegraphics[width=5cm, height=4cm]{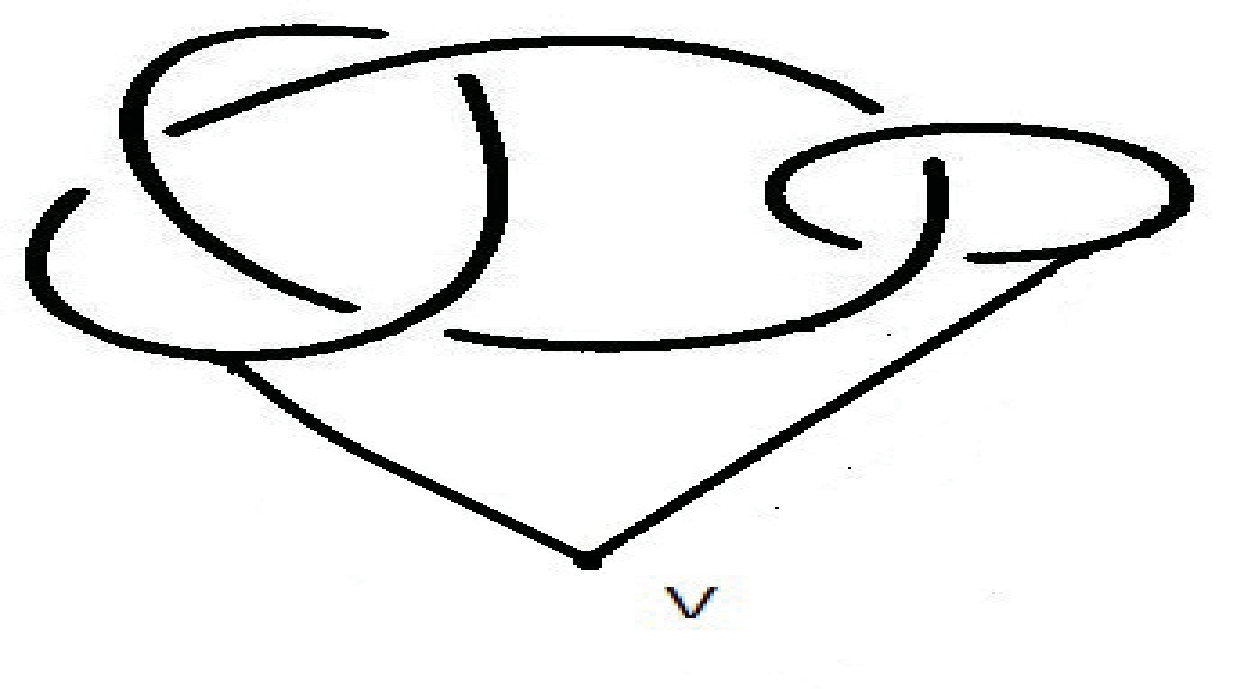}
\end{center}
\caption{A legged loop system $\gamma$ in $S^3$ 
with free fundamental group of rank $2$}
\label{fig:freeSgraph}
\end{figure}

\noindent{\bf Example~2.4.} Let $\gamma$ be a legged loop system with base point $v$ 
in $M=S^3$ illustrated in Fig.~\ref{fig:freeSgraph} with  
$\pi_1(M\setminus \gamma,v^+)$ a free group of rank $2$. In fact, 
a trivial legged loop system is obtained by sliding an edge along another edge, 
so that  $\pi_1(M\setminus \ell(\gamma),v^+)$ is a free group of rank $2$. 
A regular neighborhood $V$ of $\gamma$ in $M$ and the 
closed complement $V'=\mbox{cl}(M\setminus V)$ constitute a genus $2$ Heegaard splitting $V\cup V'$ of $M$ 
by noting that the 3-manifold $V'$ is a handlebody of genus $2$ by the loop system theorem and the Alexander theorem (cf. e.g., \cite{K1}). Thus, the union $V\cup V'$ is a genus $2$ Heegaard splitting of $M$. 
Since the legged loop system  $\gamma$ with  base point $v$ is a spine of $V$ by sliding the base point $v$ into $\partial V$, there is a Heegaard graph $\gamma_M$ 
of $M$ with $\gamma$ as $\gamma^+_V$. 
By Lemma~2.3, the spun $S^2$-link $K(\gamma_M)$ in the 4-sphere $X(M)=S^4$ has 
the free fundamental group $\pi_1(X(M)\setminus K(\gamma_M), v^+)$ of rank 4, which 
does not admit any meridian basis because the spun $S^2$-link $K(\gamma_M)$ in $S^4$ contains, as a component, the spun trefoil $S^2$-knot  whose fundamental group is known to be not infinite cyclic. 

\phantom{x}

Given a proper arc system $a_*$ in $M^{(o)}$, there is a legged loop system $\gamma$ in $M$ with 
the proper arc system $a(\gamma)=a_*$ in $M^{(o)}$. 
The spun $S^2$-link $K(\gamma)$ in $X(M)$ is uniquely determined by the arc system $a_*$ and thus denoted by $S(a_*)$. 
The following lemma is  used toward the final step of the proof of  Poinca{\'e} conjecture.

\phantom{x}

\noindent{\bf Lemma~2.5.} Let $a_*$ be a proper arc system in a compact once-punctured manifold 
$M^{(o)}=\mbox{cl}(M\setminus B)$ of a homotopy 3-sphere $M$. 
If the spun $S^2$-link $S(a_*)$ in the 4-sphere $X(M)$ is a trivial $S^2$-link, then 
the proper arc system $a_*$ is in a boundary-collar $S\times[0,1]$ of $M^{(o)}$. 

\phantom{x}

\noindent{\bf Proof of Lemma~2.5.} 
By Lemma~2.1, the fundamental group $\pi_1(M^{(o)}\setminus a(\gamma),v^+)$ is a meridian-based free group. 
Consider the 2-sphere $S$ as the boundary 
\[\partial(d\times[0,1])=d\times0\cup (\partial d)\times[0,1]\cup d\times 1\]
 of the product $d\times[0,1]$ for a disk $d$ so that 
$d\times 0$ contains one end of the proper arc system $a_*$ and $d\times 1$ contains the other end of the proper arc system $a_*$. Let $(E;E_0,E_1)$ be the triplet obtained from $(M^{(o)},d\times 0,d\times 1)$ 
by removing a tubular neighborhood of $a_*$ in $M^{(o)}$. For $v^+\in E_0$, the inclusion $E_0\subset E$ induces an 
isomorphism 
\[\pi_1(E_0,v^+)\to \pi_1(E,v^+).\]
By \cite[Theorem~10.2]{Hempel}, $E$ is diffeomorphic to the connected sum of the product $E_0\times[0,1]$ and a 
homotopy 3-sphere. 
This means that the proper arc system $a_*$ is in a boundary-collar $S\times[0,1]$. 
This completes the proof of Lemma~2.5.
$\square$

\phantom{x}

\noindent{\bf  3. A ribbon surface-link and a disk-chord system of  a ribbon handlebody system}

By combining Lemmas~2.3 with Free Ribbon Lemma in Section~1, the following lemma is obtained. 

\phantom{x}

\noindent{\bf Lemma~3.1.} 
The spun $S^2$-links $K(\gamma)$ of every Heegaard link $\gamma$ of every homotopy 3-sphere$M$ is a ribbon $S^2$ link in $X(M)$. 

\phantom{x}

The following lemma makes a connection betwen a knot in $M$ and a Heegaard graph of $M$. 

\phantom{x}

\noindent{\bf Lemma~3.2.} For every knot $k$ in every homotopy 3-sphere$M$, 
there is a Heegaard graph $\gamma$ of $M$ such that the knot $k$ is equivalent to 
a component of the loop system $\ell(\gamma)$ of $\gamma$ in  $M$. 

\phantom{x}

\noindent{\bf Proof of Lemma~3.2.} By considering $k$ as a polygonal loop in $M$, 
there is a triangulation $\mathcal T$ of $M$ whose 1-skeleton 
${\mathcal T}^{(1)}$ contains the knot $k$. The graph ${\mathcal T}^{(1)}$ is deformed 
into a legged loop system $\gamma'$ in $M$ so that $k$ is a component of 
the loop system $k(\gamma')$. 
Let $V'$ be a regular neighborhood of $\gamma'$ in $M$ which is a handlebody. 
The legged loop system $\gamma'$ is deformed into a spine $\gamma_{V'}$ of the handlebody $V'$.
The closed complement $V=\mbox{cl}(M\setminus V)$ is also a handlebody, 
so that there is a Heegaard splitting $V\cup V'$ of $M$. Hence there is  a Heegaard graph 
$\gamma$ of $M$ obtained from $\gamma_V$ and $\gamma_{V'}$ such that $k$ is equivalent to a component of the loop system $\ell(\gamma)$. 
 $\square$

\phantom{x}

By Lemma~3.2, there is a Heegaard graph 
$\gamma$ of $M$ whose loop system contains the knot $k$. 
By Lemma~3.1, the spun $S^2$-link $K(\gamma)$ is a ribbon $S^2$-link in $X(M)$, 
so that the spun $S^2$-knot $K(k)$ is a ribbon $S^2$-knot in $X(M)$ 
because any component of a ribbon $S^2$-link in $S^4$ is a ribbon $S^2$-knot 
in $S^4$ by definition.  Thus, the following result is obtained.

\phantom{x}

\noindent{\bf Lemma~3.3.} For every knot $k$ in every homotopy 3-sphere $M$, 
the spun $S^2$-knot $K(k)$  is a ribbon $S^2$-knot in $X(M)$.

\phantom{x}

For a knot $k$ in the interior of $M^{(0)}=\mbox{cl}(M\setminus B)$ for a 3-ball $B$,  
the {\it spun torus-knot} of $k$ is a torus-knot $T(k)$ in $X(M)$ given by the inclusions 
\[T(k)=k\times S^1\subset M^{(o)}\times S^1
\subset M^{(o)}\times S^1\cup S\times D^2=X(M).\] 
The spun torus-knot $T(k)$ in $X(M)$ is uniquely constructed up to choices 
of a 3-ball $B$. 
The following lemma is important to our purpose.

\phantom{x}

\noindent{\bf Lemma~3.4.} For every knot $k$ in every homotopy 3-sphere $M$, 
the spun torus-knot $T(k)$  is a ribbon torus-knot in $X(M)$. 

\phantom{x}

\begin{figure}[hbtp]
\begin{center}
\includegraphics[width=14cm, height=5cm]{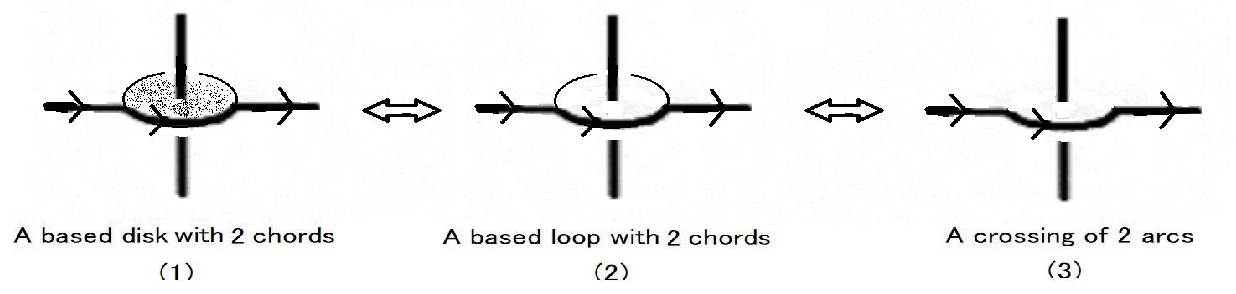}
\end{center}
\caption{Two arcs of $k$ near a disk $d_i$ drawn as thick lines}
\label{fig:2chordchange}
\end{figure}

\phantom{x}

\noindent{\bf Proof of Lemma~3.4.} 
From construction, the spun $S^2$-knot $K(k)$ in $X(M)$ is obtained from $T(k)$ by the unique 2-handle surgery, so that 
the spun torus-knot $T(k)$ is obtained from the spun $S^2$-knot $K(k)$ 
by the converse 1-handle surgery. 
By definition, the spun torus-knot $T(k)$ is a ribbon torus-knot, completing the proof.
$\square$

\phantom{x}

Assume that  a ribbon surface-link $L$ is  obtained from a trivial oriented 
$S^2$-link $O$ by surgery along a 1-handle system $h_*$ of 
disjointedly embedded oriented 1-handles $h_j\,(j=1,2,\dots,s)$ (for some $s$) 
on $O$ in $S^4$.  
A ribbon handlebody system bounded by a ribbon surface-link is discussed here 
(see\cite[II.3.61]{KSS}). 
Let $B_*$ be a system of disjoint 3-balls 
$B_i\,(i=1,2,\dots,m)$ in $S^4$ bounded by $O$. 
The intersection $h_j\cap O$ consists of two disks, called the {\it attaching disks} of 
$h_j$ to $O$. A {\it meridian disk} of the 1-handle $h_j$ is a proper disk in $h_j$ parallel to any one of the attaching disks. By an isotopic deformation of the 1-handle system  $h_*$, the intersection $h_*\cap \mbox{Int}B_i$ can be assumed to be a meridian disk system (possible empty) in $h_*$, whose number of meridian disks is called the 
{\it ribbon index} of $h_*$ in $B_i$. 
A {\it ribbon handlebody system} of a ribbon surface-link $L$ is the union 
\[ V_R=B_*\cup h_*,\] 
which is an immersed handlebody system bounded by $L$ in $S^4$. 
The {\it ribbon index} of $V_R$ is the total number of the ribbon indexes of $h_*$ 
in $B_i$ for all $i$. The {\it disk-chord system} of a ribbon surface-link $L$ is the pair 
$(d_*,\alpha_*)$ of a disk system $d_*$, called a {\it based disk system}, 
and an arc system $\alpha_*$, called a {\it chord system}, 
in $S^4$ obtained from the ribbon handlebody system $V_R=B_*\cup h_*$ 
by shrinking the 3-ball $B_i$ into a disk $d_i$ for every $i$ and then shrinking the 
1-handle $h_j$ into a core arc $\alpha_j$ of $h_j$ spanning the loop system $o_*=\partial d_*$, called a {\it based loop system}, for every $j$. 
See Fig.~\ref{fig:2chordchange} (1) for a situation around a disk in a based disk system. 
From construction, the ribbon index of $h_*$ in $B_i $ is equal to 
the number of the transverse intersection points $\alpha_*\cap \mbox{Int}d_i$, called the {\it chord index} of $\alpha_*$ in $d_i$. 
The {\it chord index} of the disk-chord system $(d_*,\alpha_*)$ is the total number of the chord indexes of $\alpha_*$ in $d_i$ for all $i$. 
By the orientations of $L$ and $S^4$, the based disk system $d_*$ can be uniquely  oriented, and the ribbon handlebody system $V_R$ and the ribbon surface-link $L$ are uniquely recovered from the disk-chord system $(d_*,\alpha_*)$ by thickening the chord system $\alpha_*$ and the based disk system $d_*$, where an argument in \cite{HK} is needed for uniqueness of the embedded 1-handle system. Let 
\[\Delta^2\subset \Delta^3\subset \Delta^4\]
be the inclusions such that $\Delta^4$ is a 4-ball in $S^4$, $\Delta^3$
is a proper 3-ball of $\Delta^4$ and  $\Delta^2$ is a proper disk of  $\Delta^3$. 
A disk-chord system $(d_*,\alpha_*)$ of $L$ in $S^4$ 
can be moved into $\mbox{Int}\Delta^3$ isotopically by first moving a neighborhood of the based disk system $d_*$ into $\mbox{Int}\Delta^3$ and then moving the remaining part of the arc system $\alpha_*$ into $\mbox{Int}\Delta^3$ (see \cite[II.3.61]{KSS}).
So, assume that a disk-chord system $(d_*,\alpha_*)$ of $L$ is in 
$\mbox{Int}\Delta^3$. 
The ribbon handlebody system $V_R$ and the ribbon surface-link $L$ are uniquely realized from a disk-chord system $(d_*,\alpha_*)$ of $L$ in $\mbox{Int}\Delta^4$. 
A {\it chord graph} of $L$ is the graph $o_*\cup \alpha_*$ in $\mbox{Int}\Delta^3$ 
obtained from a disk-chord system $(d_*,\alpha_*)$ in $\mbox{Int}\Delta^3$ 
by taking $o_*=\partial d_*$. 
A {\it chord diagram } of $L$ is a diagram $C(o_*,\alpha_*)$ in $\mbox{Int}\Delta^2$ for a chord graph $o_*\cup\alpha_*$ of $L$ in $\mbox{Int}\Delta^3$. 
A ribbon surface-link $L$ in $S^4$ is uniquely realized in 
$\mbox{Int}\Delta^4$ from a chord graph $o_*\cup \alpha_*$ of $L$ 
in $\mbox{Int}\Delta^3$ and also from a chord diagram $C(o_*,\alpha_*)$ of $L$ 
in $\mbox{Int}\Delta^2$, because the based loop system $o_*$ in 
$\mbox{Int}\Delta^3$ constructs uniquely the trivial $S^2$-link $O$ by 
the Horibe-Yanagawa lemma in \cite{KSS}. 
On the other hand, a ribbon handlebody system $V_R$ of $L$ 
cannot be uniquely recovered because in general a disjoint disk system 
$d_*$ in the interior of $\Delta^3$ 
with $\partial d_*=o_*$ is not unique (see \cite[Lemma~I.1.4]{KSS}). 
So, to fix a ribbon handlebody system $V_R$ of $L$, every loop of the based loop system $o_*$ should be fixed as it is shown in of Fig.~\ref{fig:2chordchange} (2). 
The following observation is obtained from the above argument. 

\phantom{x}

\noindent{\bf Observation~3.5.} 
A ribbon surface-link $L$ and a ribbon handlebody system $V_R$ in $S^4$ 
are uniquely realized in $\mbox{Int}\Delta^4$ from a disk-chord system $(d_*,\alpha_*)$ in 
$\mbox{Int}\Delta^3$, and also from a chord graph $o_*\cup \alpha_*$ in $\mbox{Int}\Delta^3$ or 
a chord diagram $C(o_*,\alpha_*)$ in $\mbox{Int}\Delta^2$ by fixing every loop of 
the based loop system $o_*$ as it is shown in  Fig.~\ref{fig:2chordchange} (2).

\phantom{x}

A chord diagram has the advantage of being easy to handle. For example, 
the moves on chord diagrams for equivalent ribbon surface-links are known in  \cite{K2,K3,K4, K5}. 
A ribbon handlebody $V_R$ bounded by a ribbon torus-knot $T$ is called a 
{\it ribbon solid torus}. 
The following lemma is an easy exercise of the moves on chord diagrams in \cite{K2} and used in Section~4. 

\phantom{x}

\noindent{\bf Lemma~3.6.} Every ribbon solid torus of ribbon index $n$ bounded by a ribbon torus-knot $T$ in $\mbox{Int}\Delta^4$ is deformed into a ribbon solid torus 
$V_R$ with $\partial V_R = T$ which
is realized by a disk-chord system $(d_*,\alpha_*)$  in  $\mbox{Int}\Delta^3$ of  
$\mbox{Int}\Delta^4$ where
\[d_*=\{d_i|\, i=1,2,\dots,n\},\quad 
\alpha_*=\{\alpha_i|\, i=1,2,\dots,n\}\quad \mbox{and}\quad o_*=\partial d_*\]  
such that 

\medskip 

\noindent{(1)} the chord $\alpha_i$ connects $o_i$ to $o_{i+1}$ 
for every $i\, (i=1,2,\dots,n)$ with $o_{n+1}=o_1$, and

\medskip 

\noindent{(2)} the chord index of $\alpha_*$ to $d_i$ is equal to $1$ for every $i$.

\phantom{x}

\begin{figure}[hbtp]
\begin{center}
\includegraphics[width=12cm, height=4.5cm]{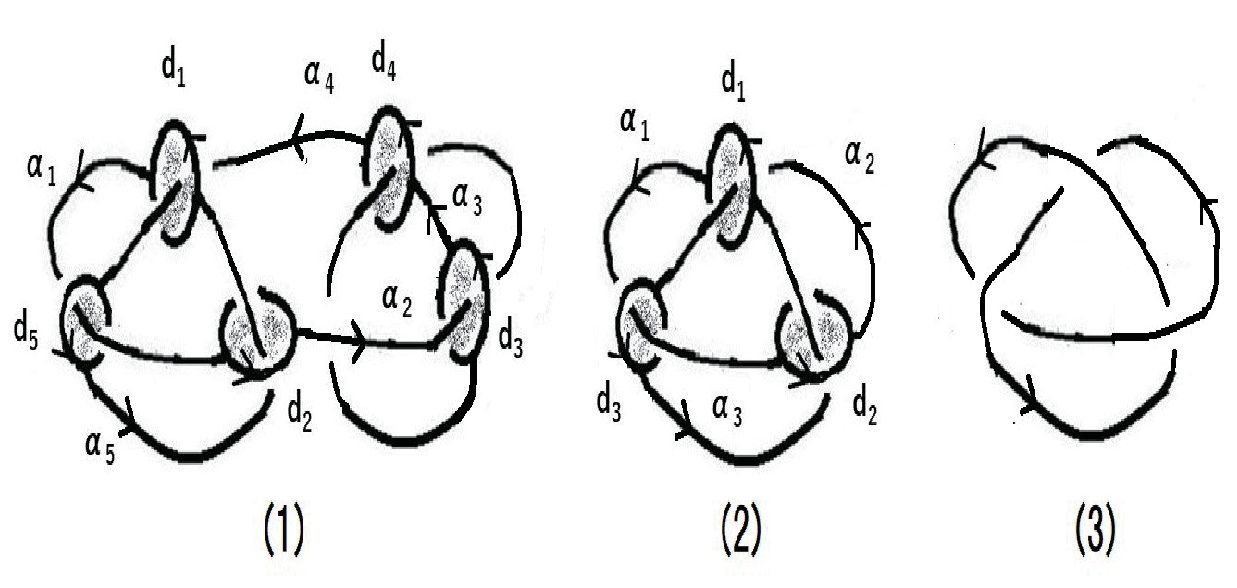}
\end{center}
\caption{CP disk-chord systems of ribbon solid tori  (1), (2) bounded 
by the spun torus-kot of the trefoil knot (3)}
\label{fig:circular}
\end{figure}

The disk-chord system $(d_*,\alpha_*)$ in Lemma~3.6 
is called  a {\it circular primitive disk-chord system} or 
briefly a {\it CP disk-chord system} (see Fig.~\ref{fig:circular} (1), (2) for examples). 
The {\it spine} of a disk-chord system $(d_*,\alpha_*)$ is a graph $\Gamma$ obtained from $d_*\cup\alpha_*$ by shrinking every disk $d_i$ into a vertex $v_i$ for every $i$. 
A {\it regular maximal tree}  of  $\Gamma$ is a  tree $\tau^+$ in $\Gamma$  
obtained from a maximal tree $\tau$ of  $\Gamma$ by taking a regular neighborhood of $\tau$ in $\Gamma$. 
A {\it regular maximal tree} of a disk-chord system $(d_*,\alpha_*)$ is a disk-chord
system $\tau^+ (d_*,\alpha_*)$ obtained from  a regular maximal tree 
$\tau^+$ of the spine $\Gamma$ by making every vertex $v_i$ in $\tau^+$ back to the original disk $d_i$ for every $i$. Let $\dot \tau^+ (d_*,\alpha_*)=\dot \tau^+$ be
the set of all the degree $1$ vertexes of $\tau^+$. 
The arc system 
\[e_*=\mbox{cl}(\Gamma\setminus \tau^+)
=\mbox{cl}((d_*\cup\alpha_*)\setminus \tau^+(d_*,\alpha_*))
\]
is called the
{\it complementary arc system} of a regular maximal tree $\tau^+(d_*,\alpha_*)$ in 
a disk-chord system  $(d_*,\alpha_*)$. 

\phantom{x}

\noindent{\bf 4. Main result: Proof of Theorem~1.1}

Throughout this section, the proof of Theorem~1.1 is done. Let $k$ be a knot in a
homotopy 3-sphere $M$. If $k$ is a trivial knot in $M$, then the knot $k$ belongs 
to a 3-ball $D^3$ in $M$. So, assume that $k$ is a non-trivial oriented knot in $M$. Since the spun torus-knot 
$T(k)$ is a ribbon torus-knot in $X(M)$ by Lemma~3.4, there is a ribbon solid torus $V_R$ of some ribbon index $n$ with 
$\partial V_R = T(k)$ in $\mbox{Int}\Delta^4$ which is realized
by a CP disk-chord system  $(d_*,\alpha_*)$ of chord index $n$ in 
$\mbox{Int}\Delta^3$  and 
a chord diagram $C(d_*,\alpha_*)$ in $\mbox{Int}\Delta^2$  by Observation~3.5. 
Since there is a meridian-preserving isomorphism  
$\pi_1(M\setminus k,v^+)\to\pi_1(X(M)\setminus T(k),v^+)$ by the
van Kampen theorem,  the longitude of $k$ in 
$M$  represents an infinite order element in the fundamental group $\pi_1(X(M)\setminus T(k),v^+)$. 
This implies that an oriented meridian loop of $V_R$ is a uniquely determined loop in 
$T(k)$ up to isotopies of $T(k)$, and the CP disk-chord system $(d_*,\alpha_*)$  is assumed that $k$ meets $d_i$ with just one boundary arc and
just one interior point transversely for every $i$, as in 
Fig.~\ref{fig:2chordchange} (1) (see also Fig.~\ref{fig:circular} (1), (2) for examples).  
Assume that $k$ is in   $\mbox{Int}M^{(o)}$.
The following lemma is obtained.

\phantom{x}

\noindent{\bf Sublemma~4.1.}
The disk system $d_i\,(i=1,2,\dots,n)$ is deformed into $\mbox{Int}M^{(o)}$ 
by an isotopy of $X(M)$ keeping the knot $k$ fixed. 

\phantom{x}

\noindent{\bf Proof of Sublemma~4.1.} 
For every $i$, let $c_i$ be a simple arc in $d_i$ connecting the point $k\cap\mbox{Int}(d_i)$ to 
a point in the arc $k\cap\partial d_i$. 
The arc system $c_i\,(i=1,2,\dots,n)$ is deformed into a bi-collar neighborhood $M^{(o)}\times [-1,1]$ of $M^{(o)}$ with 
$M^{(o)}\times 0=M^{(o)}$ in $X(M)$ by an isotopy keeping $M^{(o)}$ fixed. Then the arc system $c_i\,(i=1,2,\dots,n)$ is projected into 
$M^{(o)}$ by a general position argument. 
A deformed disk system $d_i\,(i=1,2,\dots,n)$ in $M^{(o)}$ is obtained from 
the arc system $c_i\,(i=1,2,\dots,n)$ in $M^{(o)}$ by widening them as a small disk system, completing the proof of Sublemma~4.1. $\square$

\phantom{x}

By Sublemma~4.1, consider that the CP disk-chord system 
$(d_*,\alpha_*)$ of $V_R$ is in $M^{(o)}$. 
The spine $\Gamma$ of $(d_*,\alpha_*)$ is a degree $4$ graph in $M^{(o)}$. 
For every regular maximal tree $\tau^+$ of $\Gamma$, there is a disk 
$\delta^2$ in $M^{(o)}$ with $\dot\tau^+ = \tau^+\cap  \partial \delta^2$ 
such that a neighborhood of every degree $4$ vertex of  $\tau^+$ in $\delta^2$
gives  Fig.~\ref{fig:2chordchange} (1) in $\tau^+(d_*,\alpha_*)$. 
The disk $\delta^2$ is called a {\it regular support disk} for $\tau^+(d_*,\alpha_*)$.  
This disk $\delta^2$ is moved into the 2-sphere $S = \partial M^{(o)}$. 
Let $\delta^3 = \delta^2 \times  [0, 1]$  be a collar of  $\delta^2$ in $M^{(o)}$ which is a 3-ball with $\delta^3 \cap  S = \delta^2 \times  0 = \delta^2$. 
Let $e_*$ be the complementary arc system of $\tau^+(d_*,\alpha_*)$ in 
$(d_*,\alpha_*)$ consisting of arcs $e_i\, (i = 1, 2,\dots  , n + 1)$, where $n$  is the chord index of the CP disk-chord system $(d_*,\alpha_*)$ which is  determined by the Euler characteristics $\chi(\Gamma)=-n$.  
The knot $k$ in $M^{(o)}$ is deformed in $M^{(o)}$ so that the intersection 
$t = k \cap \delta^3$ is a tangle in $\delta^3$ whose projection image  under the canonical projection
\[\delta^3 = \delta^2\times [0, 1] \to  \delta^2 \] 
is the regular maximal tree $\tau^+$ in the regular support disk $\delta^2$
by pushing $\tau^+(d_*,\alpha_*) \setminus \dot\tau^+(d_*,\alpha_*)$ into 
$\delta^2 \times  (0, 1)$ and then by creating a crossing point by the
move from (1)  to  (3) in Fig.~\ref{fig:2chordchange}.  
Then the regular maximal tree $\tau^+$ in $\delta^2$ can be regarded as a tangle diagram of $t$  in $\delta^2$. 
Let $[t, \tau^+]$ be the disk union between the tangle $t$ and the graph 
$\tau^+$ in the preimage of $\tau^+$ under the canonical projection 
$\delta^3 \to \delta^2$. 
The following sublemma is essentially observed in \cite[Theorem~2.3 (3)]{K6} for an inbound arc diagram.

\phantom{x}

\noindent{\bf Sublemma~4.2.} The spun $S^2$-link $T(t)$ of a tangle $t$ in 
$\delta^3$ in the 4-disk 
\[U^4=\delta^3\times[0,1]\times S^1\cup \delta^2\times D^2
\subset M^{(o)}\times S^1\cup S\times D^2 = X(M)\]
bounds a ribbon 3-ball system 
\[V'_R=[t,\tau^+]\times S^1\cup \tau^+\times D^2\]
which extends to a ribbon solid torus $V_R$ of the spun torus-knot $T(k)$ 
such that the compact complement $\mbox{cl}(V_R\setminus V'_R)$ 
is a disjoint 3-ball system bounded by the spun $S^2$-link $S(e_*)$ in $X(M)$.

\phantom{x}

\noindent{\bf Proof of Sublemma~4.2.} If $t$ is a 1-string tangle 
with $\tau^+$ a simple arc, then 
$V'_R=[t,\tau^+]\times S^1\cup \tau^+\times D^2$ is a 1-handle thickening $t$, that is a ribbon 3-ball with ribbon index $0$. 
If $t$ is a 2-string tangle  with $\tau^+$  just one degree $4$ vertex graph, then 
$t$ is the  2-tangle in  Fig.~\ref{fig:2chordchange} (3) and  
$V'_R$ is a ribbon 3-ball system  with ribbon index $1$ 
giving the disk chord system of  Fig.~\ref{fig:2chordchange} (1).  
In the general case of $t$ and $\tau^+$,  as  a combination result of these two observations, $V'_R$ is a ribbon 3-ball system  
giving a disk-chord system $\tau^U(d_*,\alpha_*)$ 
in the 4-disk $U^4$  such that  $ \tau^U(d_*,\alpha_*))$ is diffeomorphic to 
the regular maximal tree $\tau^+(d_*,\alpha_*)$ of $(d_*,\alpha_*)$ in $\delta^3$.  
Let $\delta^4$ be a 4-ball in $U$ with $\delta^3$ as a proper 3-ball.  
The following sublemma is needed.

\phantom{x}

\noindent{\bf Sublemma~4.3.} There is an orientation-preserving 
diffeomorphism of $X(M)$ sending
$(U^4,\tau^U(d_*,\alpha_*))$ to $(\delta^4,\tau^+(d_*,\alpha_*))$. 

\phantom{x}

\noindent{\bf Proof of Sublemma~4.3.} For the regular maximal tree $\tau^+$ in the regular support disk $\delta$, 
find a 2-disk $\delta^2_0\subset \mbox{Int}\delta$ such that 
$\tau'=\delta^2_0\cap \tau^+$ has
 $\mbox{cl}(\tau^+\setminus\tau')\cong\dot\tau^+\times[0,1]$ and construct a 4-ball 
$\delta^4_0\subset \mbox{Int}U$ with $\delta^2_0$ as a trivial proper disk.
Then construct a proper 3-ball $\delta^3_0\subset \delta^4_0$
with $\delta^2_0$ as a proper disk. 
Note that there is an orientation-preserving diffeomorphism of $S^4$ 
sending the triad $(\delta^4_0,\delta^3_0, \delta^2_0)$ to 
the triad $(\delta^4,\delta^3,\delta^2)$ and 
the regular maximal tree $\tau'(d_*,\alpha_*)$ of  $(d_*,\alpha_*)$ given by $\tau'$  in 
$\delta^3_0$ to $\tau^+(d_*,\alpha_*)$ in $\delta^3$. 
Since $\mbox{cl}(U^4\setminus \delta^4_0)$ is diffeomorphic to 
$S^3\times[0,1]$ (see \cite{K10}), there is an orientation-preserving diffeomorphism 
\[(\mbox{cl}(U^4\setminus \delta^4_0),\mbox{cl}(U^4\setminus \delta^4_0)\cap \tau^+)\to (S^3, \dot\tau^+)\times[0,1].\]
Then there is a triad $(U^4,U^3, U^2)$ with $U^3$ a proper 3-ball in $U^4$ and $U^2$ a proper 2-disk in $U^3$ such that there is an orientation-preserving diffeomorphism 
of $S^4$ sending the triad $(U^4,U^3, U^2)$  to 
the triad $(\delta^4_0,\delta^3_0,\delta^2_0)$ and $\tau^U(d_*,\alpha_*)$ 
in $U^3$ to $\tau'(d_*,\alpha_*)$ in $\delta^3_0$ . Thus, 
there is an orientation-preserving diffeomorphism of $S^4$ 
sending the triad $(U^4,U^3, U^2)$ to 
the triad $(\delta^4,\delta^3,\delta^2)$ and $\tau^U(d_*,\alpha_*)$ 
in $U^3$ to $\tau^+(d_*,\alpha_*)$ in $\delta^3$. 
This completes the proof of Sublemma~4.3. $\square$

\phantom{x}

By Sublemma~4.3, the ribbon 3-ball system 
$V'_R$ realizing $\tau^U(d_*,\alpha_*)$ in $U^4$ extends to 
a ribbon solid torus $V_R$ in $S^4$. This means that 
the spun $S^2$-link $S(e_*)$ in $X(M)$ bounds the disjoint 3-ball system
$\mbox{cl}(V_R\setminus V'_R)$. This completes the proof of Sublemma~4.2.
$\square$

\phantom{x}

By Lemma 2.5 and Sublemma~4.2, the proper arc system $e_*$ and hence $k$ are in the 3-ball $D^3$ which is a regular neighborhood 
of  $\delta^2\times [0,1]$ in $M^{(o)}$. 
This completes the proof of Theorem ~1.1. $\square$

\phantom{x}

\noindent{\bf  5. Conclusion}

A general problem arising from this paper is how 
any given ribbon solid torus bounded by the spun torus-knot $T(k)$ of a knot $k$  relates to a knot diagram $D(k)$ of $k$.  
For example, the CP disk-chord system $(d_*,\alpha_*)$  in Fig.~\ref{fig:circular} (1) is seen to represent a ribbon solid torus  bounded by the spun torus-knot $T(k)$ of the trefoil knot $k$ in Fig.~\ref{fig:circular} (3). In fact, the ribbon torus-knot given by  Fig.~\ref{fig:circular} (1)  is equivalent to the ribbon torus-knot given by Fig.~\ref{fig:circular} (2)  by  moves on chord diagrams in \cite{K2,K3,K4, K5} and 
by Sublemma~4.2 
the CP disk-chord system of Fig.~\ref{fig:circular} (2) is the CP disk-chord system of the spun ribbon solid torus of  the trefoil knot diagram $D(k)$ shown in Fig.~\ref{fig:circular} (3) .  It would be interesting to point out  that the CP disk-chord system $(d_*,\alpha_*)$ in Fig.~\ref{fig:circular} (1) is  not  
the CP disk-chord system of the spun ribbon solid torus of any  knot diagram 
$D'(k)$ of the trfoil knot $k$. 
To see this, the cross-index in \cite{KSY} is used. 
 If  $(d_*,\alpha_*)$ is obtained from the spun ribbon solid torus of a  trefoil knot diagram $D'(k)$, then the complementary arc system
$e_*$ of any regular maximal tree $\tau^+(d_*,\alpha_*)$ in  $(d_*,\alpha_*)$ in a regular support disk $\delta$ must have the  cross-index $0$ in the annulus 
$A$ given by any extended disk $\delta^+$ such that 
$\mbox{Int}\delta^+\supset \delta$ and $e$ is  an immersed arc system in the annulus  $A=(\delta^+\setminus \delta)$. However, the coss-index of $e_*$ in an annulus $A$ is $1$ for the  diagram given in Fig.~\ref{fig:circular} (1). This means that  the 
CP disk-chord system $(d_*,\alpha_*)$  in Fig.~\ref{fig:circular} (1)  is not 
the CP disk-chord system of the spun ribbon solid torus of any  trefoil knot diagram $D'(k)$. 

\phantom{x}

\noindent{\bf Acknowledgments.} For Free Ribbon Lemma, 
the author thanks the conference organizers at Sochi Conference \lq\lq Geometry and topology of 3-manifolds\rq\rq on September 18, 2022 giving motivating him to revise the proof after the author's zoom talk by hearing. 
The present version was almost completed during the author's stay at Beijing Jiaotong University China from December 2, 2023 to January 1, 2024. The author would like to 
thank Liangxia Wan and Research Assistants: Ruiyi Cui (Graduate student) and Hang Yin (Student) for providing him a quiet and comfortable stay. 
This work was partly supported by JSPS KAKENHI Grant Numbers JP19H01788, JP21H00978 and MEXT Promotion of Distinctive Joint Research Center Program JPMXP0723833165.

\phantom{x}

\end{document}